\documentclass{amsart}

\usepackage{latexsym}
\usepackage{amssymb}
\usepackage{amsmath}
\usepackage{amsthm}
\usepackage{hyperref}
\usepackage{xcolor}

\newcommand{\vvol}{\mathrm{V}}

\newcommand{\sd}{s_{\delta}}
\newcommand{\cd}{c_{\delta}}
\newcommand{\lgra}{\longrightarrow}

\newcommand{\iid}{\mathrm{Id}\,}

\newcommand{\ddiv}{\mathrm{div}}

\newcommand{\trace}{\mathrm{tr\,}}

\newcommand{\R}{\mathbb{R}}

\newtheorem{example}{Exemples}[section]
\newtheorem{thm}{Theorem}[section]

\newtheorem{lemma}[thm]{Lemma}

\newtheorem{remark}[thm]{Remark}
\newtheorem{remarks}[thm]{Remarks}
\newtheorem{definition}[thm]{Definition}
\newtheorem{notation}[thm]{Notation}
\newtheorem{exabout:ample}[thm]{Example}

\newcommand{\beqt}{\begin{equation}}  \newcommand{\eeqt}{\end{equation}}
\newcommand{\bal}{\begin{align}}      \newcommand{\eal}{\end{align}}
\newcommand{\ba}{\begin{array}}      \newcommand{\ea}{\end{array}}
\newcommand{\bc}{\begin{center}}     \newcommand{\ec}{\end{center}}
\newcommand{\be}{\begin{enumerate}}  \newcommand{\ee}{\end{enumerate}}
\newcommand{\beq}{\begin{eqnarray}}  \newcommand{\eeq}{\end{eqnarray}}
\newcommand{\beQ}{\begin{eqnarray*}} \newcommand{\eeQ}{\end{eqnarray*}}
\newcommand{\bi}{\begin{itemize}}    \newcommand{\ei}{\end{itemize}}
\newcommand{\bt}{\begin{tabular}}    \newcommand{\et}{\end{tabular}}

\begin{document}

\author[J. ROTH]{Julien Roth}
\address[J. ROTH]{Universit\'e Gustave Eiffel, CNRS, LAMA UMR 8050, F-77447 Marne-la-Vallée, France}
\email{julien.roth@univ-eiffel.fr}
\author[A. UPADHYAY]{Abhitosh Upadhyay}
 \address[A. UPADHYAY]{School of Mathematics and Computer Science, Indian Institute of Technology, Goa 403401, India}
\email{abhitosh@iitgoa.ac.in}
\thanks{The second author gratefully acknowledges the financial support from the Indian Institute of Technology Goa through Start-up Grant \textbf{2021/SG/AU/043}.}

\keywords{p-Laplacian, Steklov problem, eigenvalues}
\subjclass[2010]{53C42, 53A07, 49Q10}

\title{Reilly-type upper bounds for the $p$-Steklov problem on submanifolds }

\begin{abstract}
We prove Reilly-type upper bounds for the first non-zero eigenvalue of the Steklov problem associated with the $p$-Laplace operator on submanifolds of manifolds with sectional curvature bounded form above by a non-negative constant.
\end{abstract}
 
\maketitle

\section{Introduction}
Let $(M^n,g)$ be a $n$-dimensional compact, connected, oriented manifold without boundary, isometrically immersed into the $(n+1)$-dimensional Euclidean space $\R^{n+1}$. The spectrum of Laplacian of $(M,g)$ is an increasing sequence of real numbers
$$0=\lambda_0(\Delta)<\lambda_1(\Delta)\leqslant\lambda_2(\Delta)\leqslant\cdots\leqslant\lambda_k(\Delta)\leqslant\cdots\lgra+\infty.$$
The eigenvalue $0$ (corresponding to constant functions) is simple and $\lambda_1(\Delta)$ is the first positive eigenvalue. In \cite{Re}, Reilly proved the following well-known upper bound for $\lambda_1(\Delta)$
\begin{equation}\label{Reilly1}
\lambda_1(\Delta)\leqslant\frac{n}{\vvol(M)}\int_MH^2dv_g,
\end{equation}
where $H$ is the mean curvature of the immersion. He also proved an analogous inequality involving the higher order mean curvatures. Namely, for $r\in\{1,\cdots,n\}$
\begin{equation}\label{Reillyr}
\lambda_1(\Delta)\left(\int_MH_{r-1}dv_g\right)^2\leqslant\vvol(M)\int_MH_r^2dv_g,
\end{equation}
where $H_r$ is the $r$-th mean curvature, defined by the $r$-th symmetric polynomial of the principal curvatures.\\
\indent
Inequalities \eqref{Reilly1} and \eqref{Reillyr} have been generalized in many ways like for submanifolds of any codimension of Euclidean spaces and spheres \cite{Re,Gr}, submanifolds of hyperbolic spaces \cite{He,EI,Gr}, other differential operators of divergence-type \cite{AM,Ro2}, Paneitz-like operators \cite{Ro3} as well as for different types of Steklov problems. In particular, Ilias and Makhoul \cite{IM} proved the following upper bound for the first eigenvalue $\sigma_1$ of the Steklov problem
$$\sigma_1V(\partial M)^2\leqslant nV(M)\int_{\partial M}\|H\|^2dv_g,$$ 
where $(M^n,g)$ is a compact submanifold of $\R^N$ with boundary $\partial M$ and $H$ denote the mean curvature of $\partial M$. They also proved analogue inequalities involving higher order mean curvatures like in \eqref{Reillyr}. Recently, both the authors with Manfio have extended this inequality for submanifolds of any Riemannian manifold of bounded sectional curvature in \cite{MRU}.

Let us consider $(M^n,g)$ a compact Riemannian manifold with a possibly non-empty boundary $\partial M$. For $p\in(1,+\infty)$, we consider the so-called $p$-Laplacian defined by 
$$\Delta_pu=-\ddiv(\|\nabla u\|^{p-2}\nabla u)$$
for any $\mathcal{C}^2$ function. For $p=2$, $\Delta_2$ is nothing else than the Laplace-Beltrami operator of $(M^n,g)$.\\
Over the past years, this operator $\Delta_p$, and especially its spectrum, has been intensively studied, mainly for Euclidean domains with Dirichlet or Neumann boundary conditions (see for instance \cite{Le} and references therein) and also on Riemannian manifolds \cite{LMS}. Later, Du and Mao \cite{DM} gave analogue of Reilly inequalities \eqref{Reilly1} and \eqref{Reillyr} for the $p$-Laplacian on submanifolds of Euclidean spaces and sphere and it was extended by Chen and Wei \cite{CW} for submanifolds of the hyperbolic space. Very recently, Chen has obtained upper bounds for submanifolds of manifolde with curvature bounding from above generalizing to the $p$-Laplacian the result of Heintze for the Laplacian.\\
In the present paper, we will consider the Steklov problem associated with the $p$-Laplacian on submanifolds with boundary of the Euclidean space. Namely, we consider the $p$-Steklov problem which is the following boundary value problem

 \begin{equation}\tag{$S_p$}\label{S1}
\left\{\begin{array}{ll}
\Delta_pu=0&\text{in}\ M,\\\\
\|\nabla u\|^{p-2}\frac{\partial u}{\partial\nu}=\sigma |u|^{p-2}u&\text{on}\ \partial M,
\end{array}\right.
\end{equation}
where $\frac{\partial u}{\partial\nu}$ is the derivative of the function $u$ with respect to the outward unit normal $\nu$ to the boundary $\partial M$. Note that for $p=2$, \eqref{S1} is the usual Steklov problem (for instance one can refer to \cite{GP} regarding an overview of results about the spectral geometry of the Steklov problem). It has been observed that very little is known about the spectrum of this $p$-Steklov problem. If $M$ is a domain of $\R^N$, there exists a sequence of positive eigenvalues $\sigma_{0,p}=0<\sigma_{1,p}\leqslant\sigma_{2,p}\leqslant\cdots\leqslant \sigma_{k.p}\leqslant\cdots$ consisting in the variational spectrum and obtained by the Ljusternik-Schnirelmann theory (see \cite{Le,To} for instance). One can refer to \cite{Bo} for details about the Ljusternik-Schnirelmann principle. Note that, as mentionned in \cite[Remark 1.1]{LMS}, the arguments used in \cite{Le} can be extended to domains on Riemannian manifolds and we have that
there exists a non-decreasing sequence of variational eigenvalues obtained by the Ljusternik–Schnirelman principle. Moreover, the eigenvalue $0$ is simple with constant eigenfunctions and is isolated, that is there is no eigenvalue between $0$ and $\lambda_1$.  Then, the first positive eigenvalue of the $p$-Steklov problem $\sigma_{1,p}$ satisfies the following variational characterization
\begin{equation}\label{characsigma}
\sigma_{1,p}=\inf\left\{ \dfrac{\displaystyle\int_{M}\|\nabla u\|^pdv_g}{\displaystyle\int_{\partial M}|u|^pdv_h}\ \Bigg| u\in W^{1,p}(M)\setminus\{0\},\ \int_{\partial M}|u|^{p-2}udv_h=0\right\},
\end{equation}
where $\nabla$ is the gradient on $M$, $dv_g$ and $dv_h$ are the Riemannian volume forms respectively associated with the metric $g$ on $M$ and the induced metric $h$ on $\partial M$. \\
It is to note that all the other eigenvalues $\sigma_{k,p}$ of this sequence also have a variational characterization but we don't know if all the spectrum is contained in this sequence.\\\\
Recently,  V. Sheela  has obtained upper bound for the first eigenvalue $\sigma_{1,p}$ of the $p$-Steklov problem \eqref{S1} for Euclidean domains \cite{Sh}. Namely, she proved that for a bounded domain $\Omega$ with smooth boundary, $\sigma_{1,p}\leqslant \dfrac{1}{R^{p-1}}$ ({\it resp.} $\dfrac{n^{p-2}}{R^{p-1}}$) if $1<p<2$ ({\it resp.} $p\geqslant2$), where $R>0$ satisfies $\vvol(\Omega)=\vvol(B(R))$ and $B(R)$ is a ball of radius $R$.\\
After that, in \cite{Ro}, the first author proved upper bounds of Reilly-type for $\sigma_{1,p}$ for submanifolds with boundary of the Euclidean space. Namely, he proved that
$$\sigma_{1,p}\leqslant N^{\frac{|2-p|}{2}}n^\frac{p}{2}\left( \int_{\partial M}\|H\|^{\frac{p}{p-1}}dv_h\right)^{p-1}\frac{V(M)}{V(\partial M)^p},$$
and the more general one
$$\sigma_{1,p}\left(\int_{\partial M}\trace(T)\right)^p\leqslant N^{\frac{|2-p|}{2}}n^\frac{p}{2}\left( \int_{\partial M}\|H_T\|^{\frac{p}{p-1}}\right)^{p-1}V(M),$$
where $T$ be a symmetric and divergence-free $(1,1)$-tensor on $\partial M$. The aim of the present paper is to prove an inequality for submanifolds with boundary of Riemannian manifolds of sectional curvature bounded from above by a non-negative constant. Namely, we prove the following result.

\begin{thm}\label{thm1}
Let $\delta\geqslant0$, $p>1$ two real numbers and $(\bar{M}^N,\bar{g})$ a $N$-dimensional Riemannian manifold of sectional curvature bounded form above by $\delta$. Let $(M^n,g)$ be a compact $n$-dimensional Riemannian manifold with non empty boundary $\partial M$ isometrically immersed into $(\bar{M},\bar{g})$ and let $S$ be a symmetric, positive definite and divergence-free $(1,1)$-tensor on $\partial M$. Then $\sigma_{1,p}$ satisfies
\begin{enumerate}
\item If $\delta=0$ , then 
$$\sigma_{1,p}\left(\int_{\partial M}\trace(S)dv_h\right)^p\leqslant N^{\frac{|p-2|}{2}}n^{\frac{p}{2}}V(M)\left( \int_{\partial M}\|H_S\|^{\frac{p}{p-1}}dv_h\right)^{p-1}.$$
\item If $\delta>0$ and $M$ is contained in a ball of radius $R\leqslant \dfrac{\pi}{4\sqrt{\delta}}$, then 
\
\begin{enumerate}
\item for $1<p<2$, we have
$$ \sigma_{1,p}\leqslant\delta^{\frac{p}{2}-1}(N+1)^{\frac{2-p}{2}}n^{\frac{p}{2}}\dfrac{V(M)}{V(\partial M)}\left(\delta+\dfrac{\displaystyle\int_{\partial M}\|H_S\|^2dv_g}{ \inf(\trace(S))^2V(\partial M)}\right),$$ 
\item  and for $p\geqslant2$, we have
$$\sigma_{1,p}\leqslant(N+1)^{\frac{p-2}{2}}n^{\frac{p}{2}}\dfrac{V(M)}{V(\partial M)}\left(\delta+\dfrac{\displaystyle\int_{\partial M}\|H_S\|^2dv_g}{ \inf(\trace(S))^2V(\partial M)}\right)^{\frac{p}{2}}.$$
\end{enumerate}

\end{enumerate}
\end{thm}

\section{preliminaries}
Let $(\bar{M}^{N},\bar{g})$ be a $N$-dimensional Riemannian manifold with sectional curvature $K_{\bar{M}}\leqslant\delta$. Let $q$ a fixed point in $\bar{M}$, we denote by $r(x)$ the geodesic distance between $x$ and $q$. Moreover, we define the vector field $Z$ by 
$Z(x):=\sd(r(x))(\bar{\nabla} r)(x)$, $\sd$ is the function defined by
\[ \sd(r)=\left\{ \begin{array}{lll}

\frac{1}{\sqrt{\delta}}\sin(\sqrt{\delta}r) &\; \text{if} \;\; \delta>0\\

r &\; \text{if}\;\; \delta=0 \\

\frac{1}{\sqrt{|\delta|}}\sinh(\sqrt{|\delta|}r) &\; \text{if} \;\; \delta<0.\ \ \end{array} \right.\] 
We also define 
\[ \cd(r)=\left\{ \begin{array}{lll}

\cos(\sqrt{\delta}r) &\; \text{if} \;\; \delta>0\\

1 &\; \text{if}\;\; \delta=0 \\

\cosh(\sqrt{|\delta|}r) &\; \text{if} \;\; \delta<0.\ \ \end{array} \right.\] 
Hence, we have $\cd^2+\delta\sd^2=1$, $\sd'=\cd$ and $\cd'=-\delta\sd$.\\
In order to prove Theorem \ref{thm1}, we recall some key lemmas. First of all, we have the following Lemma, proved by Grosjean (\cite{Gr}) which in some sense gives way to extend Hsiung-Minkowki formulas in space of non constant curvature.
\begin{lemma}[\cite{Gr2}]\label{lemmaZ}
Let $(\Sigma,g)$ be a compact submanifold of $(\bar{M},\bar{g})$ and $S$ be a symmetric, positive definite and divergence-free $(1,1)$-tensor on $\Sigma$, then the following hold
\begin{enumerate}
\item $\displaystyle\sum_{i=1}^N\langle S\nabla Z_i,Z_i\rangle\leqslant \trace(S)-\delta\langle SZ^{\perp},Z^{\perp}\rangle,$
\item $\ddiv\left(SZ^{\top}\right)\geqslant n\left(\cd(r)\trace(S)+\langle Z,H_S\rangle\right).$\\\\
If in addition, $\Sigma$ has no boundary, then\\
\item $\displaystyle\int_{\Sigma}\cd(r)\trace(S)dv_g\leqslant\displaystyle\int_M\|H_S\|\sd(r)dv_g,$
\item $\delta\displaystyle\int_{\Sigma} \langle SZ^{\top},Z^{\top}\rangle dv_g\geqslant \displaystyle\int_{\Sigma}(\cd^2(r)\trace(S)-\|H_S\|\sd(r)\cd(r))dv_g$. 
\end{enumerate}
\end{lemma}
Here, $H_S$ denotes $\trace(B\circ S)$ and so is a normal vector field and $Z^{\top}$ is the part of $Z$ tangent to $\Sigma$. Note that if $S=\iid$, we recover the classical inequalities proved by Heintze \cite{He}.\\\\
To prove the desired upper bounds, we will use the variational characterization \eqref{characsigma} of $\sigma_{1,p}$. For this, we need to use appropriate test functions. As usual, for eigenvalue upper bounds for submanifolds, the candidates for test functions are the coordinates functions and their analogues in non constant curvature $Z_i=\dfrac{\sd(r)}{r}x_i$, $1\leqslant i \leqslant N$, which are the coordinates of $Z$ in a normal frame $\{e_1,\ldots,e_N\}$. To be eligible to be test functions, weed need to "center" these functions which is possible due to the following lemma given by Chen in \cite{Ch}.

\begin{lemma}[\cite{Ch}]\label{lemcenter}
Let $p\in(1,+\infty)$ and assume that $\Sigma$ is a submanifold of $\bar{M}$ contained in a convex ball $B\subset\bar{M}$. Then, there exist $q_0\in B$ such that for any $i\in\{1,\cdots,N\}$,
$$\int_{\Sigma} \left|\dfrac{\sd(r)}{r}x_i\right|^{p-2}\dfrac{\sd(r)}{r}x_idv_g=0,$$
where $r$ is the distance function to $q_0$ in $\bar{M}$.
\end{lemma}
For the case $\delta>0$, we need another test function $\cd$. In order to use it as a test function, we need to translate it to the appropriate constant. For this, we recall the following elementary lemma, also given by Chen in \cite{Ch}, for our purpose. 
\begin{lemma}\label{lemcenter2}
Let $\delta>0$, $p\in(1,+\infty)$ and assume that $\Sigma$ is a submanifold of $\bar{M}$ contained in a ball of center $q_0$ and radius $\rho<\frac{\pi}{2}$. Then, there exist a constant $c\in[0,1]$ so that
$$\displaystyle\int_{\Sigma}\left|\dfrac{\cd(r)-c}{\sqrt{\delta}}\right|^{p-2}\dfrac{\cd(r)-c}{\sqrt{\delta}}dv_g=0,$$
where $r$ is the distance function to $q_0$ in $\bar{M}$.
\end{lemma}
Finally, we recall the following technical lemma proved by Manfio and the two authors in \cite{MRU} which will be useful at the end of the proof of Theorem \ref{thm1}.
\begin{lemma}[\cite{MRU}]\label{lem2}
Let $(\bar{M}^{N},\bar{g})$ be Riemannian manifold with sectional  bounded from above by $\delta$,  $\delta>0$. Let $(\Sigma,g)$ be a closed Riemannian manifold isometrically immersed into $(\bar{M}^{N},\bar{g})$ and assume that $\Sigma$ is contained in a geodesic ball of radius $R<\frac{\pi}{2\sqrt{\delta}}$. Let $S$ be a symmetric, divergence-free and positive definite $(1,1)$-tensor on $\Sigma$. Then, we have
$$1-\left(\frac {\displaystyle\int_{\Sigma}\cd(r)dv_g}{V(\Sigma)}\right)^2\geqslant\frac{1}{1+\dfrac{\displaystyle\int_{\Sigma}\|H_S\|^2dv_g}{\delta \inf(\trace(S))^2V(\Sigma)}}.$$
\end{lemma}
\section{Proof of the Theorem \ref{thm1}}
\subsection{The case $\delta=0$}
We want to use the coordinate functions as test functions in the variational characterization of $\sigma_{1,p}$ so we need to consider the coordinates center at the good point. Therefore, we apply Lemma \ref{lemcenter} to $\Sigma=\partial M$ and we consider $r$ as the distance to the point $q_0$ given in Lemma \ref{lemcenter}.  Thus, we are able to prove the following lemma.
\begin{lemma}\label{lemsd}
For any $p\in(1,+\infty)$, we have
$$\sigma_{1,p}\int_{\partial M}r^pdv_h\leqslant N^{\frac{|p-2|}{2}}n^{\frac p2}V(M).$$
\end{lemma}
\noindent {\it Proof:} From Lemma \ref{lemcenter}, we can consider the functions $Z_i=\dfrac{\sd(r)}{r}x_i$, $1\leqslant i\leqslant N$,  as test functions in the variational characterization \eqref{characsigma} of $\sigma_{1,p}$. Since we are in the case $\delta=0$, then we have $Z_i=x_i$. Thus, taking summation for $i$ from $1$ to $N$, we have
\beq\label{varZ}
\sigma_{1,p}\int_{\partial M}\displaystyle\sum_{i=1}^N|Z_i|^pdv_h&\leqslant&\int_{M}\displaystyle\sum_{i=1}^N\|\nabla Z_i\|^pdv_g.
\eeq
Now, we will discuss the proofs for the cases $p\geqslant2$ and $1<p<2$, separately.\\\\
{\it Case $1<p<2$.}
We have 
\beqt\label{sdrp1}
r^p=\left(  r^2\right)^{\frac{p}{2}}=\left( \displaystyle\sum_{i=1}^N Z_i^2\right)^{\frac{p}{2}}\leqslant\displaystyle\sum_{i=1}^N |Z_i|^p, 
\eeqt
since $p<2$.\\
On the other hand, by the H\"older inequality (for vectors), we have
\beqt
 \displaystyle\sum_{i=1}^N \|\nabla Z_i\|^p\leqslant N^{\frac{2-p}{p}}\left(  \displaystyle\sum_{i=1}^N \|\nabla Z_i\|^2\right)^{\frac{p}{2}},
 \eeqt
which gives with the first point of Lemma \ref{lemmaZ} and $\delta=0$
\beqt\label{sdrp2}
 \displaystyle\sum_{i=1}^N \|\nabla Z_i\|^p\leqslant N^{\frac{2-p}{p}}n^{\frac{p}{2}}.
\eeqt
Hence, from \eqref{varZ}, \eqref{sdrp1} and \eqref{sdrp2}, we get
\beQ
\sigma_{1,p}\int_{\partial M}r^pdv_h&\leqslant&\sigma_{1,p}\int_{\partial M}\displaystyle\sum_{i=1}^N|Z_i|^pdv_h\\
&\leqslant&\int_{M}\displaystyle\sum_{i=1}^N\|\nabla Z_i\|^pdv_g\\
&\leqslant&N^{\frac{2-p}{2}}n^{\frac p2}V(M).
\eeQ
{\it Case $p\geqslant2$.}
By the H\"older inequality, we have
\beqt\label{sdrp3}
r^2=\displaystyle\sum_{i=1}^N| Z_i|^2\leqslant N^{\frac{p-2}{p}}\left( \displaystyle\sum_{i=1}^N| Z_i|^p\right)^{\frac{2}{p}},
\eeqt
which gives
\beqt\label{sdrp4}
r^p= N^{\frac{p-2}{2}}\left( \displaystyle\sum_{i=1}^N| Z_i|^p\right).
\eeqt
Morever, since $p\geqslant2$, we have
\beqt\label{sdrp5}
\displaystyle\sum_{i=1}^N\| \nabla Z_i\|^p\leqslant \left(\displaystyle\sum_{i=1}^N\|\nabla Z_i\|^2\right) ^{\frac{p}{2}}.
\eeqt
Finally, from \eqref{varZ}, using \eqref{sdrp4}, \eqref{sdrp5} and the first point of Lemma \ref{lemmaZ}, we get
\beQ
\sigma_{1,p}\int_{\partial M}r^pdv_h&\leqslant&\sigma_{1,p}N^{\frac{p-2}{2}}\sigma_{1,p}\int_{\partial M}\displaystyle\sum_{i=1}^N|Z_i|^pdv_h\\
&\leqslant&N^{\frac{p-2}{2}}\int_{M}\displaystyle\sum_{i=1}^N\|\nabla Z_i\|^pdv_g\\
&\leqslant&N^{\frac{p-2}{2}}\int_{M}\left(\displaystyle\sum_{i=1}^N\|\nabla Z_i\|^2\right) ^{\frac{p}{2}}dv_g\\
&\leqslant&N^{\frac{p-2}{2}}n^{\frac p2}V(M).
\eeQ
\hfill$\square$\\\\
Since $\delta=0$, then $\cd\equiv 1$ and so the third point of  Lemma \ref{lemmaZ} reduces to
$$\int_{\partial M}\trace(S)dv_h\leqslant\int_{\partial M}\sd(r)\|H_S\|dv_h.$$
Thus, we have 
\beQ
\sigma_{1,p}\left(\int_{\partial M}\trace(S)dv_h\right)^p&\leqslant& \sigma_{1,p}\left(\int_{\partial M}\sd(r)\|H_S\|dv_h\right)^p\\
&\leqslant&\sigma_{1,p}\left( \int_{\partial M}\sd^p(r)dv_h\right)\left( \int_{\partial M}\|H_S\|^{\frac{p}{p-1}}dv_h\right)^{p-1}\\
&\leqslant&N^{\frac{|p-2|}{2}}n^{\frac p2}V(M)\left( \int_{\partial M}\|H_S\|^{\frac{p}{p-1}}dv_h\right)^{p-1},
\eeQ
where we have used first the H\"older inequality and then Lemma \ref{lemsd}. 
\subsection{The case $\delta>0$}
In the case $\delta>0$, in addition to the $Z_i$'s, we need another test function. For this, from the assumption that $M$ is contained in a ball of radius $R<\frac{\pi}{4\sqrt{\delta}}$ and the Lemma \ref{lemcenter} the point $q_0$ belongs to this ball and so we can conclude that $M$ is contained in a ball of center $q_0$ and radius smaller than $\frac{\pi}{2\sqrt{\delta}}$. Therefore,  we can apply Lemma \ref{lemcenter2} to get a constant $c\in[0,1]$ so that 
$$\displaystyle\int_{\Sigma}\left|\dfrac{\cd(r)-c}{\sqrt{\delta}}\right|^{p-2}\dfrac{\cd(r)-c}{\sqrt{\delta}}dv_g=0.$$
For briefness, we set the function $C=\dfrac{\cd(r)-c}{\sqrt{\delta}}$ which can be used as a test function.  Hence, from the variational characterization \eqref{characsigma} of $\sigma_{1,p}$ using $C$ and the $Z_i$'s as test functions, we get
\beq\label{characCZ}
\sigma_{1,p}\displaystyle\int_{\partial M}\left(|C|^p+\sum_{i=1}^N|Z_i|^p\right)dv_h&\leqslant&\displaystyle\int_{M}\left(\|\nabla C\|^p+\displaystyle\sum_{i=1}^N\|\nabla Z_i\|^p\right)dv_g.
\eeq
Moreover, we have
\beq\label{CZ}
C^2+\displaystyle\sum_{i=1}^NZ_i^2&=&\left( \dfrac{\cd(r)-c}{\sqrt{\delta}}\right)^2+\displaystyle\sum_{i=1}^N\left(\frac{\sd(r)}{r}x_i\right)^2\nonumber\\
&=&\sd^2(r)+\frac{\cd^2(r)+c^2-2c\cd(r)}{\delta}\nonumber\\
&=&\frac{1+c^2-2c\cd(r)}{\delta},
\eeq
where we have used $\cd^2+\delta\sd^2=1$.\\
On the other hand, we also have
$$
\nabla C=\nabla\left( \dfrac{\cd(r)-c}{\sqrt{\delta}}\right)\\=-\sqrt{\delta}\sd(r)\nabla r=\sqrt{\delta}Z^{\perp}
$$
so that 
\beq\label{nablaCZ}
\|\nabla C\|^2+\displaystyle\sum_{i=1}^N\|\nabla Z_i\|^2&=&\delta\|Z^{\perp}\|^2+\displaystyle\sum_{i=1}^N\|\nabla Z_i\|^2\nonumber\\
&\leqslant&\delta\|Z^{\perp}\|^2+(n-\delta\|Z^{\perp}\|^2)\nonumber\\
&=&n,
\eeq
where we have used the first point of Lemma \ref{lemmaZ}. Note that here, $Z^{\top}$ is the part of $Z$ tangent to $M$.\\
From now on, we will consider the cases $1<p<2$ and $p\geqslant 2$ separately.\\
{\it Case $1<p<2$.} Since $p<2$, we have
\beq\label{sumCZ1}
|C|^p+\displaystyle\sum_{i=1}^N |Z_i|^p&=&\frac{1}{\delta^{\frac{p}{2}}}\left( \left|\cos(\sqrt{\delta}\,r)-c\right|^p+\displaystyle\sum_{i=1}^N \left|\sin(\sqrt{\delta}\,)\frac{x_i}{r}\right|^p\right).
\eeq
Since $\left|\sin(\sqrt{\delta}\,)\frac{x_i}{r}\right|\leqslant 1$, $\left|\cos(\sqrt{\delta}\,r)-c\right|<1$ and $1<p<2$, we have 
$$\left|\sin(\sqrt{\delta}\,)\frac{x_i}{r}\right|^p\geqslant \left|\sin(\sqrt{\delta}\,)\frac{x_i}{r}\right|^2\ \text{and}\left|\cos(\sqrt{\delta}\,r)-c\right|^p\geqslant \left|\cos(\sqrt{\delta}\,r)-c\right|^2,$$
which after reporting into \eqref{sumCZ1} gives
\beq\label{sumCZ1bis}
|C|^p+\displaystyle\sum_{i=1}^N |Z_i|^p&\geqslant&\frac{1}{\delta^{\frac{p}{2}}}\left( \left|\cos(\sqrt{\delta}\,r)-c\right|^2+\displaystyle\sum_{i=1}^N \left|\sin(\sqrt{\delta}\,)\frac{x_i}{r}\right|^2\right)\nonumber\\
&=&\frac{1}{\delta^{\frac{p}{2}-1}}\left( C^2+\displaystyle\sum_{i=1}^N Z_i^2\right)\nonumber\\
&=&\frac{1}{\delta^{\frac{p}{2}}}\Big(1+c^2-2c\cd(r)\Big),
\eeq where we have used \eqref{CZ} for the last line.
 On the other hand, by the H\"older inequality, we have
\beq\label{sumnablaCZ1}
\|\nabla C\|^p+\displaystyle\sum_{i=1}^N \|\nabla Z_i\|^p&\leqslant&(N+1)^{\frac{2-p}{2}}\left( \|\nabla C\|^2+\displaystyle\sum_{i=1}^N \|Z_i\|^2\right)^{\frac{p}{2}}\nonumber\\
&\leqslant&(N+1)^{\frac{2-p}{2}}n^{\frac{p}{2}}
\eeq
by using \eqref{nablaCZ}. Thus, using \eqref{characCZ} together with \eqref{sumCZ1bis} and \eqref{sumnablaCZ1}, we get
\beq\label{majsigma1}
\sigma_{1,p}\displaystyle\int_{\partial M}\Big(1+c^2-2c\cd(r)\Big)dv_h&\leqslant&\delta^{\frac{p}{2}}(N+1)^{\frac{2-p}{2}}n^{\frac{p}{2}}V(M).
\eeq
Moreover we have
\beq\label{ineqintcd}
\int_{\partial M}\Big(1+c^2-2c\cd(r)\Big)dv_h&=&V(\partial M)\left(1+c^2-2c\frac{\int_{\partial M}\cd(r)dv_h}{V(\partial M)}\right)\nonumber\\
&=&V(\partial M)\left(1+\left(c-\frac{\int_{\partial M}\cd(r)dv_h}{V(\partial M)}\right)^2-\left(\frac{\int_{\partial M}\cd(r)dv_h}{V(\partial M)}\right)^2\right)\nonumber\\
&\geqslant&V(\partial M)\left(1-\left(\frac{\int_{\partial M}\cd(r)dv_h}{V(\partial M)}\right)^2\right).
\eeq
Reporting this into \eqref{majsigma1}, we get
\beq\label{majsigma11}
\sigma_{1,p}V(\partial M)\left(1-\left(\frac{\int_{\partial M}\cd(r)dv_h}{V(\partial M)}\right)^2\right)dv_h&\leqslant&\delta^{\frac{p}{2}}(N+1)^{\frac{2-p}{2}}n^{\frac{p}{2}}V(M).
\eeq
Finally, we conclude by applying Lemma \ref{lem2} to get
\beq\label{majsigma12}
\sigma_{1,p}&\leqslant&\delta^{\frac{p}{2}-1}(N+1)^{\frac{2-p}{2}}n^{\frac{p}{2}}\frac{V(M)}{V(\partial M)}\left(\delta+\dfrac{\displaystyle\int_{\partial M}\|H_S\|^2dv_g}{ \inf(\trace(S))^2V(\partial M)}\right).
\eeq
{\it Case $p\geqslant2$.} On one hand, since $p\geqslant 2$, we have
\beq\label{sumnablaCZ2}
\|\nabla C\|^p+\displaystyle\sum_{i=1}^N \|Z_i\|^p&\leqslant&\left(\|\nabla C\|^2+\displaystyle\sum_{i=1}^N \|Z_i\|^2\right)^{\frac{p}{2}}\nonumber\\
&\leqslant&n^{\frac{p}{2}}
\eeq
where we have used \eqref{nablaCZ}. On the other hand, we have by the H\"older inequality
\beq\label{sumCZ2}
\left(C^2+\displaystyle\sum_{i=1}^N Z_i^2\right)^{\frac{p}{2}}&\leqslant&(N+1)^{\frac{p-2}{2}}\left(|C|^p+\displaystyle\sum_{i=1}^N |Z_i|^p\right).
\eeq
Thus, using successively \eqref{CZ}, \eqref{sumCZ2}, \eqref{characCZ} and \eqref{sumnablaCZ2}, we get
\beq\label{majsigma13}
\sigma_{1,p}\displaystyle\int_{\partial M}\left( \frac{1+c^2-2c\cd(r)}{\delta}\right)^{\frac{p}{2}}dv_h&=&\sigma_{1,p}\displaystyle\int_{\partial M}\left( C^2+\displaystyle\sum_{i=1}^N Z_i^2\right)^{\frac{p}{2}}dv_h\nonumber\\
&\leqslant&\sigma_{1,p}(N+1)^{\frac{p-2}{2}}\displaystyle\int_{\partial M}\left( |C|^p+\displaystyle\sum_{i=1}^N |Z_i|^p\right)dv_h\nonumber\\
&\leqslant&(N+1)^{\frac{p-2}{2}}\displaystyle\int_{M}\left( \|\nabla C\|^p+\displaystyle\sum_{i=1}^N \|\nabla Z_i\|^p\right)dv_h\nonumber\\
&\leqslant&(N+1)^{\frac{p-2}{2}}n^{\frac{p}{2}}V(M).
\eeq
In addition, from the H\"older inequality (for integrals), we have
\beq\label{holdercd}
\displaystyle\int_{\partial M}\left( \frac{1+c^2-2c\cd(r)}{\delta}\right)^{\frac{p}{2}}dv_h\nonumber \\ \geqslant V(\partial M)^{\frac{2-p}{2}}\left(\displaystyle\int_{\partial M}\frac{1+c^2-2c\cd(r)}{\delta}dv_h\right)^{\frac{p}{2}}.
\eeq
Hence, we deduce from \eqref{majsigma13} with \eqref{ineqintcd} and \eqref{holdercd}
\beq\label{majsigma14}
\sigma_{1,p}\dfrac{V(\partial M)}{\delta^{\frac{p}{2}}}\left[ 1-\left(\frac{\int_{\partial M}\cd(r)dv_h}{V(\partial M)}\right)^2\right]^{\frac{p}{2}}&\leqslant&\sigma_{1,p}V(\partial M)^{\frac{2-p}{p}}\left(\displaystyle\int_{\partial M}\frac{1+c^2-2c\cd(r)}{\delta}dv_h\right)^{\frac{p}{2}}\nonumber\\
&\leqslant&(N+1)^{\frac{p-2}{2}}n^{\frac{p}{2}}V(M).
\eeq
Finally, we use Lemma \ref{lem2} to conclude that
$$\sigma_{1,p}\leqslant (N+1)^{\frac{p-2}{2}}n^{\frac{p}{2}}\dfrac{V(M)}{V(\partial M)}\left(\delta+\dfrac{\displaystyle\int_{\partial M}\|H_S\|^2dv_g}{ \displaystyle\inf_{\partial M}(\trace(S))^2V(\partial M)}\right)^{\frac{p}{2}}.$$
This concludes the proof of Theorem \ref{thm1}. \hfill$\square$

\section{New results for the $p$-Laplacian when $\delta>0$}
We finish by giving similar results fo the first eigenvalue of the $p$-Laplacian for closed submanifolds when $\delta>0$. We will not give all the details of the proof since it is similar to the proof of Theorem \ref{thm1}. The difference is that the variational characterization of $\lambda_{1,p}$ is the following
\begin{equation}\label{characlambda}
\lambda_1=\inf\left\{ \dfrac{\displaystyle\int_{M}\|\nabla u\|^pdv_g}{\displaystyle\int_{M}|u|^pdv_g}\ \Bigg| u\in W^{1,p}(M)\setminus\{0\},\ \int_{M}|u|^{p-2}udv_g=0\right\}.
\end{equation}
In this case $M$ has no boundary and so we can apply Lemma \ref{lemcenter} and Lemma \ref{lemcenter2} with $\Sigma=M$ to be able to use the functions $Z_i$, $1\leqslant i\leqslant N$ and $C$ as test functions. By completely similar computations, we obtain the analogue of inequality \eqref{majsigma11} if $1<p<2$, that is,
$$\lambda_{1,p}V(M)\left(1-\left(\frac{\int_{M}\cd(r)dv_h}{V(M)}\right)^2\right)dv_h\leqslant\delta^{\frac{p}{2}}(N+1)^{\frac{2-p}{2}}n^{\frac{p}{2}}V(M),$$
and of \eqref{majsigma14} if $p\leqslant2$, that is,
$$\lambda_{1,p}\dfrac{V(M)}{\delta^{\frac{p}{2}}}\left[ 1-\left(\frac{\int_{M}\cd(r)dv_h}{V(M)}\right)^2\right]^{\frac{p}{2}}\leqslant(N+1)^{\frac{p-2}{2}}n^{\frac{p}{2}}V(M).$$
Finally, applying Lemma \ref{lem2} to $M$, we deduce the following result
\begin{thm}\label{thm2}
Let $\delta>0$, $p\in(1,+\infty)$ and $(\bar{M}^N,\bar{g})$ a Riemannian manifold of sectional curvature bounded form above by $\delta$. Let $(M^n,g)$ be a closed Riemannian manifold  isometrically immersed into $(\bar{M},\bar{g})$ and $S$ a symmetric, positive definite and divergence-free $(1,1)$-tensor on $\partial M$. We denote by $\lambda_{1,p}$ the first eigenvalue of the $p$-Laplacian on $M$. If $M$ is contained in a ball of radius $R\leqslant \dfrac{\pi}{4\sqrt{\delta}}$ then we have
\begin{enumerate}
\item If $1<p<2$, 
$$\lambda_{1,p}\leqslant\delta^{\frac{p}{2}-1}(N+1)^{\frac{2-p}{2}}n^{\frac{p}{2}}\left(\delta+\dfrac{\displaystyle\int_{M}\|H_S\|^2dv_g}{ \displaystyle\inf_M(\trace(S))^2V(M)}\right).$$
\item If $p\geqslant2$,
$$\lambda_{1,p}\leqslant(N+1)^{\frac{p-2}{2}}n^{\frac{p}{2}}\left(\delta+\dfrac{\displaystyle\int_{M}\|H_S\|^2dv_g}{ \displaystyle\inf_M(\trace(S))^2V(M)}\right)^{\frac{p}{2}}.$$
\end{enumerate}
\end{thm}
If $S=\iid$, we recover the result of Chen \cite{Ch}.

\end{document}